\documentclass[leqno]{amsart}

\usepackage{amsmath,amsfonts,amsthm,amssymb, mathrsfs, mathabx, accents, mathdots, vntex}

\usepackage[french,english]{babel}
\usepackage[bookmarks]{hyperref}
\hypersetup{backref, colorlinks=true}

 \usepackage{graphicx}
\usepackage[all]{xy}
\newtheorem{theorem}{Theorem}[section]

\newtheorem{proposition}[theorem]{Proposition}
\newtheorem{corollary}[theorem]{Corollary}
\newtheorem{remark}[theorem]{Remark}

\newtheorem{example}[theorem]{Example}
\newtheorem{definition}[theorem]{Definition}

\newtheorem{conjecture}[theorem]{Conjecture}
\newtheorem{question}[theorem]{Question}

\newenvironment{preuve}{{\em{\noindent \textbf{Proof.} }}}
{\hfill $\blacksquare$}

\def\C{ \mathbb{C}}

\def\R{ \mathbb{R}}

\def\Rank{ {\rm Rank}}

\def\rond{\mathaccent"7017}

\begin{document}

\large

\title[]{\small A REMARK ON A POLYNOMIAL MAPPING FROM $\C^n$ TO $\C^{n-1}$}
\makeatother

\author[Nguyen Thi Bich Thuy]{Nguyen Thi Bich Thuy}
\address[Nguyen Thi Bich Thuy]{UNESP, Universidade Estadual Paulista, ``J\'ulio de Mesquita Filho'', S\~ao Jos\'e do Rio Preto, Brasil}
\email{bichthuy@ibilce.unesp.br}

\maketitle \thispagestyle{empty}

\begin{abstract} 
We provide relations of the results obtained  in the articles \cite{ThuyCidinha} and  \cite{VuiThang}. Moreover, we provides some examples to illustrate these relations, using the software {\it Maple} to complete the complicate calculations of the examples. 
We give some discussions on these relations. 

\end{abstract}

\section{INTRODUCTION}

In \cite{ThuyCidinha}, we construct singular varieties ${\mathcal{V}}_G$ associated to a polynomial mapping $G : \C^{n} \to \C^{n - 1}$ where $n \geq 2$ such that if $G$ is a local submersion but is not a fibration, then the 2-dimensional   homology and intersection homology (with total perversity) 
of the variety ${\mathcal{V}}_G$ are not trivial. 
In \cite{VuiThang}, the authors  prove that if there exists a so-called {\it very good projection with respect to the regular value $t^0$} of a polynomial mapping $G: \C^n \to \C^{n-1}$, then this value is an atypical value of $G$ if and only if the Euler characteristic of the fibers is not constant. 
This paper provides relations of the results obtained  in the articles \cite{ThuyCidinha} and  \cite{VuiThang}. Moreover, we provides some examples to illustrate these relations, using the software {\it Maple} to complete the calculations of the examples.
We provide some discussions on these relations. 

\section{PRELIMINARIES}

\subsection{Intersection homology.}
We briefly recall the definition of intersection homology; for details,
we refer  to the fundamental work of M. Goresky and R. MacPherson
\cite{GM1} (see also \cite{Jean Paul}).

\begin{definition} 
{\rm Let $X$ be a $m$-dimensional variety.  A {\it stratification of $X$} is the data of a finite  filtration 
$$X = X_{m} \supset X_{m-1} \supset \cdots \supset X_0 \supset X_{-1} = \emptyset,$$
 such that  for every $i$,  the set $S_i = X_i\setminus X_{i-1}$ is either an emptyset or a manifold of dimension $i$. A connected component of $S_i$ is  called   {\it a stratum} of $X$.}
\end{definition}

We denote by $cL$  the open cone on the space $L$, the cone on the empty set being a point. Observe that if $L$ is a stratified set then $cL$ is stratified by  the cones over the strata of $L$ and an additional  $0$-dimensional stratum (the vertex of the cone). 

\begin{definition}
{\rm A stratification of $X$ is said to be {\it locally topologically trivial} if for every $x \in X_i\setminus X_{i-1}$, $i \ge 0$, there is an open neighborhood $U_x$  of $x$ in $X$, a stratified set $L$ and a homeomorphism 
 $$h:U_x \to (0;1)^i \times  cL,$$   such  that
 $h$ maps the strata of $U_x$ (induced stratification) onto the strata of  $  (0;1)^i \times cL$  (product stratification).}
\end{definition}

The definition of perversities has originally been given by Goresky and MacPherson:
\begin{definition}
{\rm 
A {\it perversity} is an $(m+1)$-uple  of integers $\bar p = (p_0, p_1, p_2,
p_3,\dots , p_m)$ such that $p_0 = p_1 = p_2 = 0$ and $p_{k+1}\in\{p_k, p_k + 1\}$, for $k \geq 2$.

Traditionally we denote the zero perversity by
$\overline{0}=(0, 0, \dots,0)$, the maximal perversity by
$\overline{t}=(0, 0, 0,1,\dots,m-2)$, and the middle perversities by
$\overline{m}=(0, 0, 0,0,1,1,\dots, [\frac{m-2}{2}])$ (lower middle) and
$\overline{n}=(0, 0, 0,1,1,2,2,\dots ,[\frac{m-1}{2}])$ (upper middle). We say that the
perversities $\overline{p}$ and $\overline{q}$ are {\it complementary} if $\overline{p}+\overline{q}=\overline{t}$.

Let $X$ be a  variety such that $X$ admits a  locally topologically trivial stratification. We say that an $i$-dimensional subset $Y\subset X$ is  $(\bar
p, i)$-{\it allowable} if  
$$\dim (Y \cap X_{m-k}) \leq i - k + p_k \text{ for
all } k.$$

Define $IC_i ^{\overline{p}}(X)$ to be the $\R$-vector subspace of $C_i(X)$
consisting in the chains $\xi$ such that $|\xi|$ is
$(\overline{p}, i)$-allowable and $|\partial \xi|$ is
$(\overline{p}, i - 1)$-allowable.}
\end{definition}

\begin{definition} 
{\rm The {\it $i^{th}$ intersection homology group with perversity $\overline{p}$}, denoted by
$IH_i ^{\overline{p}}(X)$, is the $i^{th}$ homology group of the
chain complex $IC^{\overline{p}}_*(X).$
}
\end{definition}

The notation $IH_*^{\overline{p}, c}(X)$ will refer to the intersection homology with compact supports, the notation $IH_*^ {\overline{p},cl}(X)$  will refer to the intersection homology with closed supports. In the compact case, they coincide and will be denoted by $IH_*^{\overline{p}}(X)$. 
In general, when we write $H_*(X)$ ({\it resp.,} $IH_*^{\overline{p}}(X)$), we mean the homology ({\it resp.}, the intersection homology) with both compact supports and closed supports.

 Goresky and MacPherson proved that the intersection homology  is independent on the
choice of the stratification satisfying the locally topologically trivial conditions \cite{GM1, GM2}. 

The Poincar\'e duality holds for the intersection homology of a (singular) variety:

 \begin{theorem}[Goresky,
MacPherson \cite{GM1}] {\it For any orientable compact stratified
semi-algebraic $m$-dimensional variety  $X$,  the generalized Poincar\'e duality holds:
$$IH_k ^{\overline{p}}(X)  \simeq IH_{m-k} ^{\overline{q}} (X),$$
where $\overline{p}$ and $\overline{q}$ are complementary perversities.
}
\end{theorem}

\noindent For the non-compact case, we have:
$$IH_k ^{\overline{p}, c}(X)  \simeq IH_{m-k} ^{\overline{q}, cl} (X).$$

\subsection{The bifurcation set, the set of asymptotic critical values and the asymptotic set} 
Let $G : \C^n \to \C^m$ where $n \geq m$ be a polynomial mapping. 

\medskip

i) The bifurcation set of $G$, denoted by $B(G)$ is the smallest set in $\C^m$ such that $G$ is not $C^{\infty}$ - fibration on this set (see, for example, \cite{Ku}). 

\medskip 


ii) When $n = m$, we denote by $S_G$ the set of points at which the mapping $G$ is not proper, {\it i.e.} 
$$S_G:= \{ \alpha \in \C^m : \exists \{ z_k\} \subset \C^n, \vert z_k \vert \to \infty \text{ such that } G(z_k) \to \alpha\},$$
and call it the {\it asymptotic variety} (see \cite{Jelonek}). The following holds: $B(G) = S_G$ (\cite{Jelonek}). 


\section{Varieties  ${\mathcal{V}}_G$ associated to a polynomial mapping $G: \C^n \to \C^{n-1}$} \label{VG}

In \cite{ThuyCidinha}, we construct singular varieties associated to a polynomial mapping $G: \C^n \to \C^{n-1}$ as follows:  let $G: \C^n \to \C^{n-1}$ such that $K_0(G) = \emptyset$, where $K_0(G)$ is  the set of critical values of $G$. 
Let  $\rho: \C^n \to \R$ be a real function such that 
$$\rho = a_1 \vert z_1 \vert ^2 + \cdots + a_n \vert z_n \vert^2,$$ 
where $\sum_{i = 1}^n  a^2_i \neq 0$, $\, a_i \geq 0$ 
and $a_i \in \R$.
 Let us denote $\varphi = \frac{1}{1 + \rho}$ and consider $(G, \varphi)$  as a real mapping from $\R^{2n}$ to $\R^{2n-1}$. 
 Let us define
$${\mathcal{M}}_G : = Sing(G, \varphi) = \{x \in \R^{2n} \text{ such that } \Rank_{\R} D(G, \varphi)(x) \leq 2n - 2\},$$
where $D(G, \varphi)(x)$ is the (real) Jacobian matrix of $(G, \varphi): \R^{2n} \to \R^{2n - 1}$ at $x$. 
Notice that $ Sing(G, \varphi) = Sing(G, \rho)$, so we have ${\mathcal{M}}_G=  Sing(G, \rho)$.

\begin{proposition} \cite{ThuyCidinha} \label{prodimSingGpsi}
{\rm For an open and dense set of polynomial mappings $G: \C^n \to \C^{n-1}$ such that $K_0(G) = \emptyset$, the variety ${\mathcal{M}}_G$  is a smooth manifold of  dimension $2n-2$. 
}
\end{proposition}

 Now, let us consider:

a) $F : = G_{\vert {\mathcal{M}}_G }$ the restriction of $G$ on ${\mathcal{M}}_G$,

\medskip 

b) ${\mathcal{N}}_G = {\mathcal{M}}_G \setminus F^{-1}(K_0(F))$.

\medskip

\noindent  Since the dimension of ${\mathcal{M}}_G$  is $2n - 2$ (Proposition \ref{prodimSingGpsi}), then locally, in a neighbourhood of any point $x_0$ in ${\mathcal{M}}_G$, 
   we get a mapping $F : \R^{2n-2}  \to \R^{2n - 2}$. Then 
there exists a co\-ver\-ing $\{ U_1, \ldots , U_p \}$ of ${\mathcal{N}}_G$  by open semi-algebraic subsets (in $\R^{2n}$) such that on every element of this co\-ver\-ing, the mapping $F$ induces a diffeomorphism   
onto its image (see Lemma 2.1 of \cite{Valette}). We can find semi-algebraic closed  subsets $V_i \subset U_i$ (in ${\mathcal{N}}_G$) which cover ${\mathcal{N}}_G$ as well. 
Thanks to Mostowski's Separation Lemma (see Separation Lemma in \cite{Mos}, page 246), for each $\, i =1, \ldots , p$, 
there exists a Nash function $\psi_i : {\mathcal{N}}_G \to \R$,  
such that  $\psi_i$ is positive on $V_i$ and negative on ${\mathcal{N}}_G \setminus U_i$. 
 We can choose the Nash functions $\psi_i$ such that $\psi_i (x_k)$ tends to  zero when $\{x_k\} \subset {\mathcal{N}}_G$ tends to infinity.
Let the Nash functions $\psi_i$ and $\rho$ be such that $\psi_i (x_k)$ tends to  zero and $\rho(x_k)$ tends to infinity  when $x_k \subset {\mathcal{N}}_G$ tends to infinity. 
Define a variety ${\mathcal{V}}_G$ associated to $(G, \rho)$ as 
$${\mathcal{V}}_G : = \overline{(F, \psi_1, \ldots, \psi_p)({\mathcal{N}}_G)},$$
that means  ${\mathcal{V}}_G$ is the closure of ${\mathcal{N}}_G$ by $(F, \psi_1, \ldots, \psi_p)$. 

In order to understand better the construction of the variety ${\mathcal{V}}_G$, see the example 4.13 in \cite{ThuyCidinha}.

\begin{proposition} \cite{ThuyCidinha} \label{remarkSG}
{\rm Let $G: \C^n \to \C^{n-1}$ be a polynomial mapping such that $K_0(G) = \emptyset$ and let $\rho: \C^n \to \R$ be a real  function such that 
$$\rho = a_1 \vert z_1 \vert ^2 + \cdots + a_n \vert z_n \vert^2,$$
where $\sum_{i = 1}^n a^2_i \neq 0$, $\, a_i  \geq 0$ and $a_i \in \R$ for $i = 1, \ldots, n.$ Then, there exists a real algebraic variety ${\mathcal{V}}_G$ in $\R^{2n - 2 + p}$, where $p > 0$, such that: 
 
1) The real dimension of ${\mathcal{V}}_G$ is $2n -2$,

2) The singular set at infinity 
of the variety ${\mathcal{V}}_G$  is contained in ${\mathcal{S}}_G(\rho) \times \{0_{\R^{ p}}\},$ 
where 
$${\mathcal{S}}_G(\rho):=\{ \alpha \in \C^{n-1} \, \vert \, \exists \{z_k\} \subset Sing(G, \rho) :  z_k \text{ tends to infinity}, G(z_k) \text{ tends to } \alpha \}.$$
}
\end{proposition}

\section {The Bifurcation set $B(G)$ and the homology, intersection homology of varieties  ${\mathcal{V}}_G$ associated to a polynomial mapping $G: \C^n \to \C^{n-1}$} 

We have the two following theorems dealing with the homology and intersection ho\-mo\-logy of the variety ${\mathcal{V}}_G$.

\begin{theorem} \cite{ThuyCidinha} \label{ThuyCidinha1} 
{\rm Let $G = (G_1, G_2): \C^3 \rightarrow \C^2$ be a polynomial mapping such that $K_0(G) = \emptyset$. If $B(G) \neq \emptyset$ then

\begin{enumerate}
\item[1)] $ H_2({\mathcal{V}}_G, \R) \ne 0,$
\item[2)] $IH_2^{\overline{t}}({\mathcal{V}}_G, \R) \ne 0,$ 
 where $\overline{t}$ is the total perversity.
\end{enumerate}
}
\end{theorem}



\begin{theorem} \cite{ThuyCidinha} \label{ThuyCidinha2} 
{\rm 
Let $G = (G_1, \ldots, G_{n-1}): \C^n \to \C^{n-1}$, where $n \geq 4$, be a polynomial mapping such that $K_0(G) = \emptyset$ and $\Rank_{\C} {(D \hat {G_i})}_{i=1,\ldots,n-1} \geq n-2$, where $\hat {G_i}$ is the leading form of $G_i$, that is the homogenous part of highest degree of $G_i$,  for $i = 1, \ldots, n-1$. If $B(G) \neq \emptyset$ then 

\begin{enumerate}
\item[1)] $ H_2({\mathcal{V}}_G, \R) \ne 0,$
\item[2)] $H_{2n-4}({\mathcal{V}}_G, \R) \ne 0,$

\item[3)] $IH_2^{\overline{t}}({\mathcal{V}}_G, \R) \ne 0,$ 
 where $\overline{t}$ is the total perversity.
\end{enumerate}
}
\end{theorem} 

\begin{remark} 
{\rm The singular set at infinity of ${\mathcal{V}}_G$ depends on the choice of the function $\rho$, since when  $\rho$ changes, the set ${\mathcal{S}}_G$ also changes. However, we have alway the property $B(G) \subset {\mathcal{S}}_G(\rho)$  (see \cite{Cidinha}).
}
\end{remark}

\begin{remark}
{\rm The variety ${\mathcal{V}}_G$ depends on the choice of the function $\rho$ and the functions $\psi_i$, but the theorems \ref{ThuyCidinha1}  and \ref{ThuyCidinha2} do not depend on the varieties ${\mathcal{V}}_G$. Form now, we denote by ${\mathcal{V}}_G(\rho)$ any variety ${\mathcal{V}}_G$ associated to $(G, \rho)$. If we refer to ${\mathcal{V}}_G$, that means a variety ${\mathcal{V}}_G$ associated to $(G, \rho)$ for any $\rho$.  
}
\end{remark}

\section{The Bifurcation set $B(G)$ and the Euler characteristic  of the fibers of a polynomial mapping $G: \C^n \to \C^{n-1}$} \label{relationTC-VT}

Let $G = (G_1, G_2, \ldots, G_{n-1}) : \C^n \to \C^{n-1}$ be a non-constant polynomial mapping and $t^0 = (t^0_1, t^0_2, \ldots, t^0_{n-1}) \in \C^{n-1}$ be a regular value of $G$.

\begin{definition} \cite{VuiThang}
{\rm A linear function $L :  \C^n \to  \C$ 
is said to be a very good projection with respect to the value $t^0$ 
if there exists a positive number $\delta$ 
such that for all 
$t \in D_\delta (t^0)= \{ t=(t_1,t_2, \ldots, t_{n-1} ) \in \C^{n-1} : |t_i -t_i^0| < \delta \}$:

(i) The restriction $L_t := {L |}_{ G^{-1}(t)}  \to \C$ is proper,

(ii) The cardinal  of $ L^{-1}(\lambda)$ does not depend on $\lambda$, where $\lambda$ is a regular value of $L$.}
\end{definition}

\begin{theorem} \cite{VuiThang} \label{thVuiThang}
{\rm Let $t^0$ be a regular value of $G$. Assume that there exists a very good projection with respect to the value $t^0$. Then, $t^0$ is an atypical value of $G$ if and only if the Euler characteristic of $G^{-1}(t^0)$
 is bigger than that of the generic fiber.
}
\end{theorem}

\begin{theorem} \cite{VuiThang} \label{conditionprojection}
{\rm Assume that the zero set 
$\{ z \in \C^n : \hat{G}_i (z) =0, i = 1, \ldots, n-1 \}$, where $\hat{G}_i$ is the leading form of $G_i$, 
has complex dimension one. Then any generic linear mapping $L$
 is  a very good projection with respect to any regular value $t^0$ of $G$.   
}
\end{theorem}

\section{Relations between \cite{ThuyCidinha} and \cite{VuiThang}} 
 \label{ThuyCidinhaVuiThang} 
Let $G = (G_1, \ldots, G_{n-1}): \C^n \to \C^{n-1}$ $(n \geq 3)$ be a polynomial mapping such that $K_0(G) = \emptyset$. Then any $t^0 \in \C^{n-1}$ is a regular value of $G$. Let $\rho: \C^n \to \R$ be a real function such that $\rho = a_1 \vert z_1 \vert ^2 + \cdots + a_n \vert z_n \vert^2,$ 
where $\sum_{i = 1}^n a^2_i \neq 0$, $\, a_i  \geq 0$ and $a_i \in \R$ for $i = 1, \ldots, n.$ 
 From theorems \ref{ThuyCidinha1} and \ref{thVuiThang}, we have the following corollary.  

\begin{corollary} \label{coThuyCidinha1}
{\rm Let $G = (G_1, G_2): \C^3 \rightarrow \C^2$ be a polynomial mapping such that $K_0(G) = \emptyset$. 
Assume that there exists a very good projection 
with respect to $t^0 \in \C^2$. 
If the Euler characteristic of $G^{-1}(t^0)$
 is bigger than that of the generic fiber, then 
\begin{enumerate}
\item[1)] $ H_2({\mathcal{V}}_G(\rho), \R) \ne 0,$ for any $\rho$,
\item[2)] $IH_2^{\overline{t}}({\mathcal{V}}_G(\rho), \R) \ne 0,$  for any $\rho$,
 where $\overline{t}$ is the total perversity.
\end{enumerate}
}
\end{corollary}

\begin{preuve}
Let $G = (G_1, G_2): \C^3 \rightarrow \C^2$ be a polynomial mapping such that $K_0(G) = \emptyset$. Then every point $t^0 \in \C^2$ is a regular point of $G$. 
Assume that there exists a very good projection 
with respect to $t^0 \in \C^2$. If the Euler characteristic of $G^{-1}(t^0)$
 is bigger than that of the generic fiber, then by the theorem \ref{thVuiThang}, the bifurcation set $B(G)$ is not empty. Then by the theorem \ref{ThuyCidinha1}, we have 
$ H_2({\mathcal{V}}_G(\rho), \R) \ne 0,$ for any $\rho$ and 
$IH_2^{\overline{t}}({\mathcal{V}}_G(\rho), \R) \ne 0,$  for any $\rho$,
 where $\overline{t}$ is the total perversity.
\end{preuve}

From theorems \ref{ThuyCidinha2} and \ref{thVuiThang}, we have the following corollary. 

\begin{corollary} \label{coThuyCidinha2}
{\rm Let $G = (G_1, \ldots, G_{n-1}): \C^n \to \C^{n-1}$, where $n \geq 4$, be a polynomial mapping such that $K_0(G) = \emptyset$ and $\Rank_{\C} {(D \hat {G_i})}_{i=1,\ldots,n-1} \geq n-2$, where $\hat {G_i}$ is the leading form of $G_i$. 
Assume that there exists a very good projection 
with respect to $t^0 \in \C^{n-1}$.  
 If the Euler characteristic of $G^{-1}(t^0)$
 is bigger than that of the generic fiber, then 
\begin{enumerate}
\item[1)] $ H_2({\mathcal{V}}_G(\rho), \R) \ne 0,$ for any $\rho$,
\item[2)] $H_{2n-4}({\mathcal{V}}_G(\rho), \R) \ne 0$, for any $\rho$, 
\item[3)] $IH_2^{\overline{t}}({\mathcal{V}}_G, \R) \ne 0,$  for any $\rho$, 
 where $\overline{t}$ is the total perversity.
\end{enumerate}
}
\end{corollary}
\begin{preuve}
Let $G = (G_1, \ldots, G_{n-1}): \C^n \to \C^{n-1}$, where $n \geq 4$, be a polynomial mapping such that $K_0(G) = \emptyset$. Then every point $t^0 \in \C^{n-1}$ is a regular point of $G$. Assume that there exists a very good projection 
with respect to $t^0 \in \C^{n-1}$. By the theorem \ref{thVuiThang}, the bifurcation set $B(G)$ is not empty. If $\Rank_{\C} {(D \hat {G_i})}_{i=1,\ldots,n-1} \geq n-2$, then by the theorem \ref{ThuyCidinha2}, we have 
\begin{enumerate}
\item[1)] $ H_2({\mathcal{V}}_G(\rho), \R) \ne 0,$ for any $\rho$,
\item[2)] $H_{2n-4}({\mathcal{V}}_G(\rho), \R) \ne 0$, for any $\rho$, 
\item[3)] $IH_2^{\overline{t}}({\mathcal{V}}_G, \R) \ne 0,$  for any $\rho$, 
 where $\overline{t}$ is the total perversity.
\end{enumerate}
\end{preuve}

We have also the following corollary. 

\begin{corollary} \label{coThuyCidinha3}
{\rm Let $G = (G_1, \ldots, G_{n-1}): \C^n \to \C^{n-1}$, where $n \geq 4$, be a polynomial mapping such that $K_0(G) = \emptyset$. Assume that the zero set 
$\{ z \in \C^n : \hat{G}_i (z) =0, i = 1, \ldots, n-1 \}$ has complex dimension one, where $\hat{G}_i$ is the leading form of $G_i$. 
 If the Euler characteristic of $G^{-1}(t^0)$
 is bigger than that of the generic fiber, where $t^0 \in \C^{n-1}$, then 

\begin{enumerate}
\item[1)] $ H_2({\mathcal{V}}_G(\rho), \R) \ne 0,$ for any $\rho$,
\item[2)] $H_{2n-4}({\mathcal{V}}_G(\rho), \R) \ne 0$, for any $\rho$,
\item[3)] $IH_2^{\overline{t}}({\mathcal{V}}_G(\rho), \R) \ne 0,$ for any $\rho$, 
 where $\overline{t}$ is the total perversity.
\end{enumerate}
 }
\end{corollary}

\begin{preuve} 
At first, since  the zero set 
$\{ z \in \C^n : \hat{G}_i (z) =0, i = 1, \ldots, n-1 \}$ has complex dimension one, then by the theorem  \ref{conditionprojection}, any generic linear mapping $L$
 is  a very good projection with respect to any regular value $t^0$ of $G$. Moreover, 
 the complex dimension of the set $\{ z \in \C^n : \hat{G}_i (z) =0, i = 1, \ldots, n-1 \}$ is the complex {\it corank} of ${(D \hat {G_i})}_{i=1,\ldots,n-1}$. Then $\Rank_{\C} {(D \hat {G_i})}_{i=1,\ldots,n-1} = n-2$. By the corollary \ref{coThuyCidinha2}, we get the proof of the corollary \ref{coThuyCidinha3}. 
\end{preuve}

\begin{remark}
{\rm We can construct the variety ${\mathcal{V}}_G(L)$, where $L$ is a very good projection defined in  \ref{thVuiThang} as the following: Let $G = (G_1, \ldots, G_{n-1}): \C^n \to \C^{n-1}$, where $n \geq 2$, be a polynomial mapping such that $K_0(G) = \emptyset$. Assume that  there exists a very good projection $L: \C^n \to \C$  
with respect to $t^0 \in \C^{n-1}$. Then $L$ is a linear function. Assume that $L = \sum_{i = 1}^n a_i z_i$. Then the variety ${\mathcal{V}}_G(L)$ is defined as the variety ${\mathcal{V}}_G(\rho)$, where 
 $$\rho = \sum_{i = 1}^n \vert a_i \vert \vert z_i \vert,$$ 
with $\vert a_i \vert$, $\vert z_i \vert$ are the modules of the complex numbers $a_i$ and $z_i$, respectively. With this variety ${\mathcal{V}}_G(L)$, all the results in the corollaries \ref{coThuyCidinha1}, \ref{coThuyCidinha2} and \ref{coThuyCidinha3} hold. Moreover,  the varieties  ${\mathcal{V}}_G(L)$ makes the corollaries \ref{coThuyCidinha1}, \ref{coThuyCidinha2} and \ref{coThuyCidinha3} simplier. 

}
\end{remark}

\begin{remark} \label{remark1}
{\rm In the construction of the variety ${\mathcal{V}}_G$ \cite{ThuyCidinha} (see section \ref{VG}), if  we replace $F$ by the restriction of $(G, \varphi)$ to ${\mathcal{M}}_G$, that means

 $$F : = (G, \varphi)_{\vert {\mathcal{M}}_G },$$
then we have the same results than in \cite{ThuyCidinha}. 
 In fact, in this case, since the dimension of ${\mathcal{M}}_G$  is $2n - 2$, then locally, in a neighbourhood of any point $x_0$ in ${\mathcal{M}}_G$, 
   we get a mapping $F : \R^{2n-2}  \to \R^{2n - 1}$. There exists also a co\-ver\-ing $\{ U_1, \ldots , U_p \}$ of ${\mathcal{N}}_G$  by open semi-algebraic subsets (in $\R^{2n}$) such that on every element of this co\-ver\-ing, the mapping $F$ induces a diffeomorphism   
onto its image. We can find semi-algebraic closed  subsets $V_i \subset U_i$ (in ${\mathcal{N}}_G$) which cover ${\mathcal{N}}_G$ as well. 
Thanks to Mostowski's Separation Lemma, for each $\, i =1, \ldots , p$, 
there exists a Nash function $\psi_i : {\mathcal{N}}_G \to \R$,  
such that  $\psi_i$ is positive on $V_i$ and negative on ${\mathcal{N}}_G \setminus U_i$.  
Let the Nash functions $\psi_i$ and $\rho$ be such that $\psi_i (z_k)$ and $\varphi(z_k) = \frac{1}{1 + \rho(z_k)}$  tend to  zero  where $\{z_k\}$ is a sequence in ${\mathcal{N}}_G$ tending to infinity. 
Define a variety ${\mathcal{V}}_G$ associated to $(G, \rho)$ as 
$${\mathcal{V}}_G : = \overline{(F, \psi_1, \ldots, \psi_p)({\mathcal{N}}_G)} = \overline{(G, \varphi, \psi_1, \ldots, \psi_p)({\mathcal{M}}_G)}.$$ 
We get the $(2n-2)$-dimensional singular variety ${\mathcal{V}}_G$ in $\R^{2n - 1 +p}$, the singular set at infinity of which is  ${\mathcal{S}}_G \times \{ 0_{\R^{p + 1}} \}$.

With this construction of the set ${\mathcal{V}}_G$, the corrolaries \ref{coThuyCidinha1}, \ref{coThuyCidinha2} and \ref{coThuyCidinha3} also hold. 
}
\end{remark}

\section{Some discussions} 

A natural question is to know if the converses of the corollaries \ref{coThuyCidinha1} and  \ref{coThuyCidinha2} hold. That means, let $G = (G_1, \ldots, G_{n-1}): \C^n \to \C^{n-1}$ $(n \geq 3)$ be a polynomial mapping such that $K_0(G) = \emptyset$ then 

\begin{question}
{\rm If there exists a very good projection with respect to  $t^0 \in \C^{n-1}$ and if  either  $IH_2^{\overline{t}}({\mathcal{V}}_G, \R) \neq 0$ or  $H_2({\mathcal{V}}_G, \R) \neq 0$, then is the Euler characteristic of $G^{-1}(t^0)$
  bigger than the one of the generic fiber?
}
\end{question}

By the theorem \ref{thVuiThang},
the above question is equivalent to the following question: 

\begin{question} \label{question2}
{\rm If $B(G) = \emptyset$ then are $IH_2^{\overline{t}}({\mathcal{V}}_G, \R) = 0$ and $IH_2^{\overline{t}}({\mathcal{V}}_G, \R) = 0$?}
\end{question}
This question is equivalent to the converse of the theorems \ref{ThuyCidinha1} and  \ref{ThuyCidinha2}. Note that by the proposition \ref{remarkSG}, the singular set at infinity 
of the variety ${\mathcal{V}}_G$  is contained in ${\mathcal{S}}_G(\rho) \times \{0_{\R^{1 + p}}\}.$ 
 Moreover, in the proofs of the theorems \ref{ThuyCidinha1} and  \ref{ThuyCidinha2}, we see that the characteristics of the homology and intersection homology of the variety ${\mathcal{V}}_G(\rho)$ depend on the set ${\mathcal{S}}_G(\rho)$. 
In \cite{ThuyCidinha}, we provided an example to show that the answer to  the question \ref{question2} is negative. In fact, let 
$$G: \C^3 \to \C^2, \quad G(z, w, \zeta) = (z, z\zeta^2 + w),$$
then $K_0(G) = \emptyset$ and $B(G) = \emptyset.$ 
 if we choose the function $\rho = \vert \zeta \vert ^2,$
then  
${\mathcal{S}}_G(\rho) = \emptyset $ and  $IH_2^{\overline{t}}({\mathcal{V}}_G(\rho), \R) = 0$; 
 if we choose the function 
$\rho' = \vert w \vert ^2,$
then  ${\mathcal{S}}_G(\rho') \neq \emptyset$  and  $IH_2^{\overline{t}}({\mathcal{V}}_G(\rho'), \R) \neq 0.$  
 Then, we suggest the two following conjectures. 
\begin{conjecture}
{\rm Does there exist a function $\rho$ such that if $B(G) = \emptyset$ then ${\mathcal{S}}_G(\rho) = \emptyset$? }
\end{conjecture}

\begin{conjecture}
{\rm 
 Let $G = (G_1, \ldots, G_{n-1}): \C^n \to \C^{n-1}$ $(n \geq 2)$ be a polynomial mapping such that $K_0(G) = \emptyset$. Assume that there exists a very good projection 
with respect to $t^0 \in \C^{n-1}$.  
 If the Euler characteristic of $G^{-1}(t)$ is constant, for any $t \in \C^{n-1}$,  then there exists a real positive function $\rho: \C^n \to \R$ such that $ H_2({\mathcal{V}}_G(\rho), \R) = 0$ and  $IH_2^{\overline{t}}({\mathcal{V}}_G(\rho), \R) = 0$.
}

\end{conjecture}



\begin{remark} \label{remark2}
{\rm 
The construction of the variety ${\mathcal{V}}_G$ in \cite{ThuyCidinha} (see section \ref{VG}) can be applied for any polynomial mapping $G: \C^n \to \C^{m}$, where $1\leq m \leq n-2$, such that $K_0(G) = \emptyset$. In fact, if $G$ is generic then similarly to the propositon \ref{prodimSingGpsi}, the variety 
$${\mathcal{M}}_G : = Sing(G, \varphi) = \{x \in \R^{2n} \text{ such that } \Rank D_{\R}(G, \varphi)(x) \leq 2m \},$$
has the real dimension $2m$. 
Hence, if we consider 
$F : = G_{\vert {\mathcal{M}}_G }$, that means $F$ is  the restriction of $G$ to ${\mathcal{M}}_G$, then locally
   we get a real mapping $F : \R^{2m}  \to \R^{2m}$. 
Moreover, in this case, we also have  $ B(G) \subset {\mathcal{S}}_G(\rho)$ for any $\rho$ (see \cite{Cidinha}), where 
$${\mathcal{S}}_G(\rho):=\{ \alpha \in \C^{m} \, \vert \, \exists \{z_k\} \subset Sing(G, \rho) :  z_k \text{ tends to infinity}, G(z_k) \text{ tends to } \alpha \}.$$
So, we can use the same arguments in \cite{ThuyCidinha}, and we have the following results. 

\begin{proposition} 
{\rm Let $G: \C^n \to \C^{m}$ be a polynomial mapping, where $1 \leq m \leq n-2$, such that $K_0(G) = \emptyset$. Let $\rho: \C^n \to \R$ be a real function such that 
$$\rho = a_1 \vert z_1 \vert ^2 + \cdots + a_n \vert z_n \vert^2,$$
where $\sum_{i = 1}^n a^2_i \neq 0$, $\, a_i  \geq 0$ and $a_i \in \R$ for $i = 1, \ldots, n.$ Then, there exists a real variety ${\mathcal{V}}_G$ in $\R^{2m + p}$, where $p > 0$, such that: 
 
1) The real dimension of ${\mathcal{V}}_G$ is $2m$,

2) The singular set at infinity 
of the variety ${\mathcal{V}}_G$  is contained in ${\mathcal{S}}_G(\rho) \times \{0_{\R^{p}}\}.$ 
}
\end{proposition}

Similarly to \cite{ThuyCidinha}, we have the two following theorems (see theorems \ref{ThuyCidinha1}  and \ref{ThuyCidinha2}). 

\begin{theorem} 
{\rm 
Let $G = (G_1, G_2): \C^n \rightarrow \C^2$, where $n \geq 4$, be a polynomial mapping such that $K_0(G) = \emptyset$. If $B(G) \neq \emptyset$ then

\begin{enumerate}
\item[1)] $ H_2({\mathcal{V}}_G, \R) \ne 0,$
\item[2)] $IH_2^{\overline{t}}({\mathcal{V}}_G, \R) \ne 0,$ 
 where $\overline{t}$ is the total perversity.
\end{enumerate}
}
\end{theorem}

\begin{theorem} 
{\rm Let $G = (G_1, \ldots, G_m): \C^n \to \C^{m}$, where $3 \leq m \leq n-2$, 
be a polynomial mapping such that $K_0(G) = \emptyset$. 
Assume that $\Rank_{\C} {(D \hat {G_i})}_{i=1,\ldots,m} \geq m-1$, where $\hat {G_i}$ is the leading form of $G_i$.  
 If $B(G) \neq \emptyset$ then 
\begin{enumerate}
\item[1)] $ H_2({\mathcal{V}}_G(\rho), \R) \ne 0,$ for any $\rho$,
\item[2)] $H_{2n-4}({\mathcal{V}}_G(\rho), \R) \ne 0$, for any $\rho$,
\item[3)] $IH_2^{\overline{t}}({\mathcal{V}}_G(\rho), \R) \ne 0,$ for any $\rho$, 
 where $\overline{t}$ is the total perversity.
\end{enumerate}
}
\end{theorem} 

}
\end{remark}


\section{Examples}
\begin{example}  \label{ex}

{\rm 
We give here an example to illustrate the calculations of the set  ${\mathcal{V}}_G$
 in the case of a polynomial mapping  $G: \C^2 \to \C$
 where $K_0(G) = \emptyset$, $B(G) \neq \emptyset$ and there exists a very good projection with respect to any point of $B(G)$. In general, the calculations of the set ${\mathcal{V}}_G$  are enough complicate, but the software {\it  Maple} may support us. That is what we do in the this example.

Let us consider the Broughton's example \cite{Broughton}: 
$$G: \C^2 \to \C, \quad \quad \quad  G(z, w) = z + z^2w.$$
 We have $K_0(G) = \emptyset$ and $B(G) \neq \emptyset$. In fact, since the system of equations
$\frac{\delta G}{\delta z} = \frac{\delta G}{\delta w} = 0$
has no solutions, then $K_0(G) = \emptyset$.  Moreover, 
$$G^{-1}(0) = \{ (z, w) : z = 0 \text{ or } zw = -1 \}  \cong \C \sqcup (\C \setminus \{ 0 \}),$$
and for any $\epsilon \neq 0$, we have
$$G^{-1}(\epsilon) = \{ (z, w) : z \neq 0 \text{ and } w = (\epsilon - z) / z^2 \}  \cong  \C \setminus \{ 0 \}.$$
So $G^{-1}(0)$ is not homeomorphic to $G^{-1}(\epsilon)$ for any $\epsilon \neq 0$. Hence $B(G) = \{ 0 \}$.  We determine now all the possible {\it very good projections} of $G$ with respect to $t^0 = 0 \in B(G)$. 
In fact, for any $\delta > 0$ 
and for any $t \in D_\delta(0)$, 
we have 
$$G^{-1}(t) = \{ (z, w) \in \C^2: z + z^2w = t \neq 0 \} 
= \left\{ (z, w) \in \C^2 : z \neq 0 \text{ and }  w = \frac{t - z}{z^2} \right\}.$$
Assume that $\{(z_k, w_k)\}$ is a sequence in $G^{-1}(t)$ tending to infinity. If $z_k$ tends to infinity then $w_k$ tends to zero. If $w_k$ tends to infinity then $z_k$ tends to zero. If $L$ is a very good projection with respect to $t^0 = 0$ then, by definition, the restriction $L_t := {L |}_{ G^{-1}(t)}  \to \C$ is proper. Then $L = az + bw$, where $a \neq 0$ and $b \neq 0$. We check now the cardinal   $\sharp L^{-1}(\lambda)$  of $L^{-1}(\lambda)$ where $\lambda$ is a regular value of $L$. 
 Let us replace $w = \frac{t - z}{z^2}$ in the equation $az + bw = \lambda$, we have the following equation 
$$az + b \frac{t - z}{z^2} = 0,$$
where $z \neq 0$. This equation has always three (complex) solutions. Thus, the number  $\sharp L^{-1}(\lambda)$ does not depend on $\lambda$. Hence, any linear  function of the form $L = az + bw$, where $a \neq 0$ and $b \neq 0$, is a very good projection of $G$ with respect to $t^0 = 0$. It is easy to see that the set of very good projections of $G$ with respect to $t^0 = 0$ is dense in the set of linear functions.

We choose $L = z + w$ and we compute the variety ${\mathcal{V}}_G$ associated to $(G, \rho)$ where $\rho = \vert z \vert^2 + \vert w \vert^2$. Let us denote 
$$z = x_1 + ix_2, \quad \quad \quad w=x_3 + ix_4,$$
where $x_1, x_2, x_3, x_4 \in \R$. Consider $G$ as a real polynomial mapping, we have 
$$G(x_1, x_2, x_3, x_4) = (x_1 + x_1^2 x_3 - x_2^2x_3 - 2x_1x_2x_4, x_2 + 2x_1x_2x_3+x_1^2x_4 - x_2^2x_4),$$
and 
$$\varphi = \frac{1}{1 + \rho} = \frac{1}{ \vert z \vert^2 + \vert w \vert^2} = \frac{1}{x_1^2 + x_2^2 + x_3^2 + x_4^2}.$$

The set  ${\mathcal{N}}_G = Sing(G, \rho)$ is the set of the solutions of the determinant of the minors $3 \times 3$ of the matrix
$$D_{\R}(G, \rho) ={ \begin{pmatrix} 1 + 2x_1x_3 - 2 x_2x_4 & -2x_2x_3 - 2x_1x_4 & x_1^2 - x_2^2 & -2x_1x_2  \\  2x_2x_3 + 2x_1x_4 & 1 + 2x_1x_3 - 2x_2x_4 & 2x_1x_2 & x_1^2 - x_2^2 \\ x_1 & x_2 & x_3 & x_4\end{pmatrix}}.$$

Using {\it Maple}, we: 

A) Calculate the determinants of the minors $3 \times 3$ of the matrix  $D_{\R}(G, \rho) $.

1) Calculate the determinant of the minor defined by the columns 1, 2 and 3:

\begin{figure}[h!]
\includegraphics[scale=0.7]{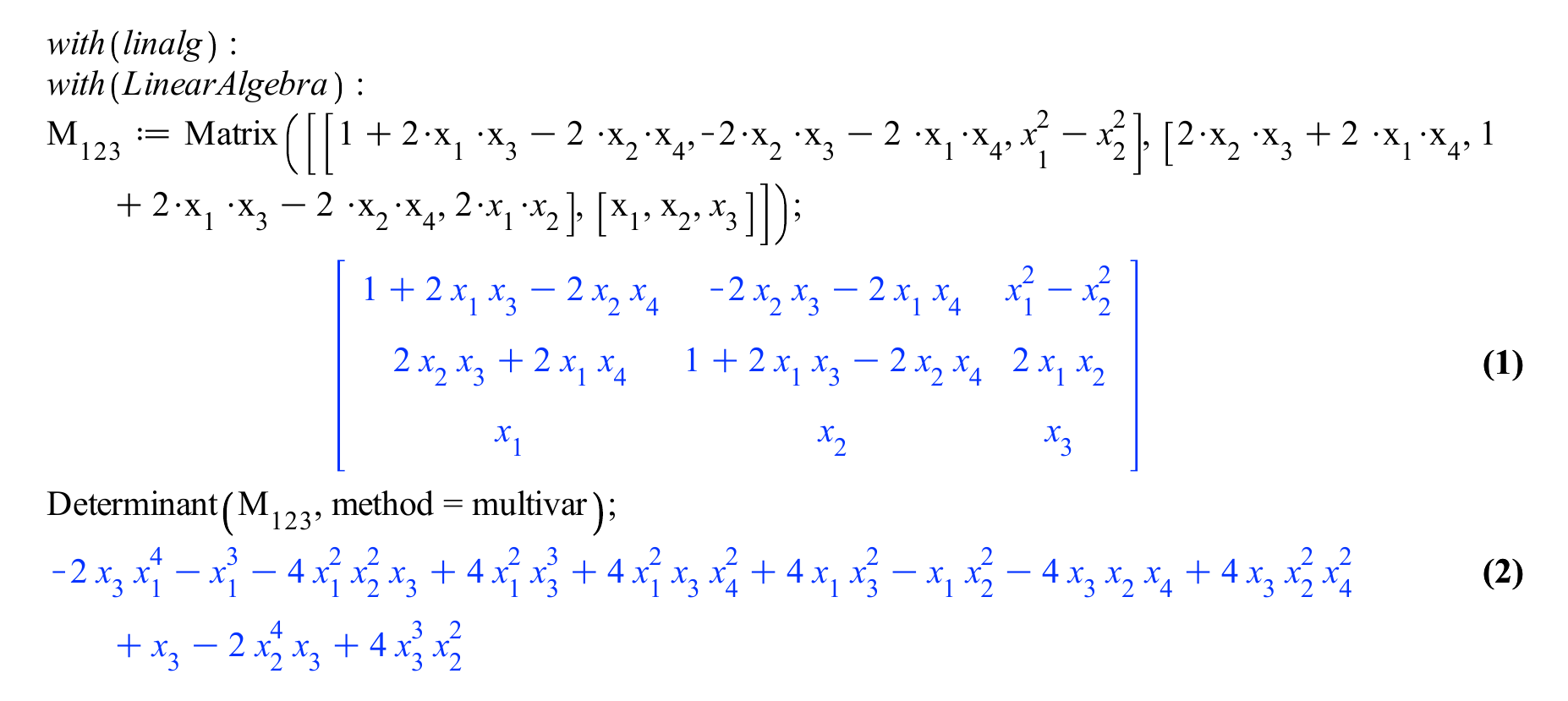}
\end{figure}

\pagebreak
2) Calculate the determinant of the minor defined by the columns 1, 2 and 4:

\begin{figure}[h!]
\includegraphics[scale=0.6]{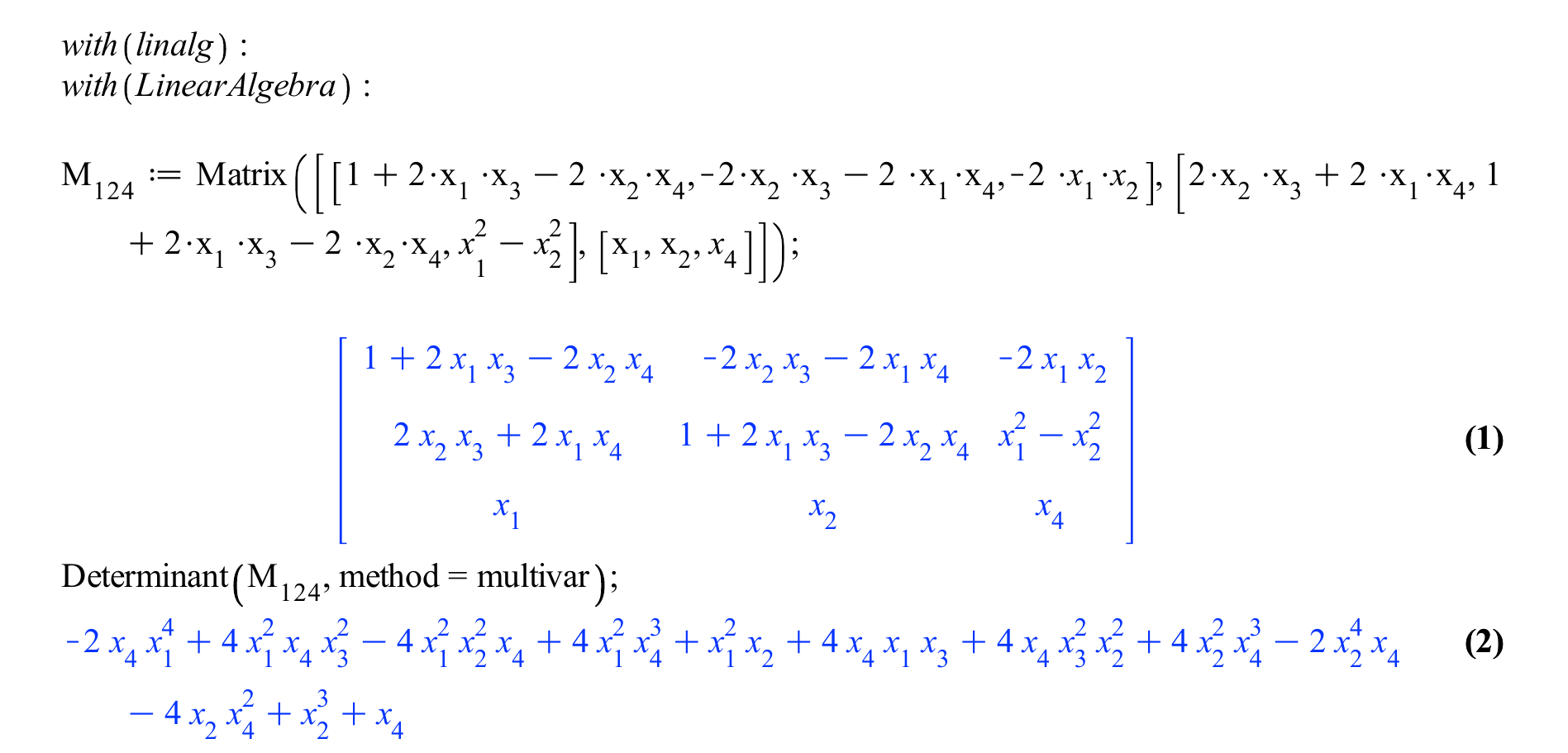}
\end{figure}

3) Calculate the determinant of the minor defined by the columns 1, 3 and 4:

\begin{figure}[h!]
\includegraphics[scale=0.6]{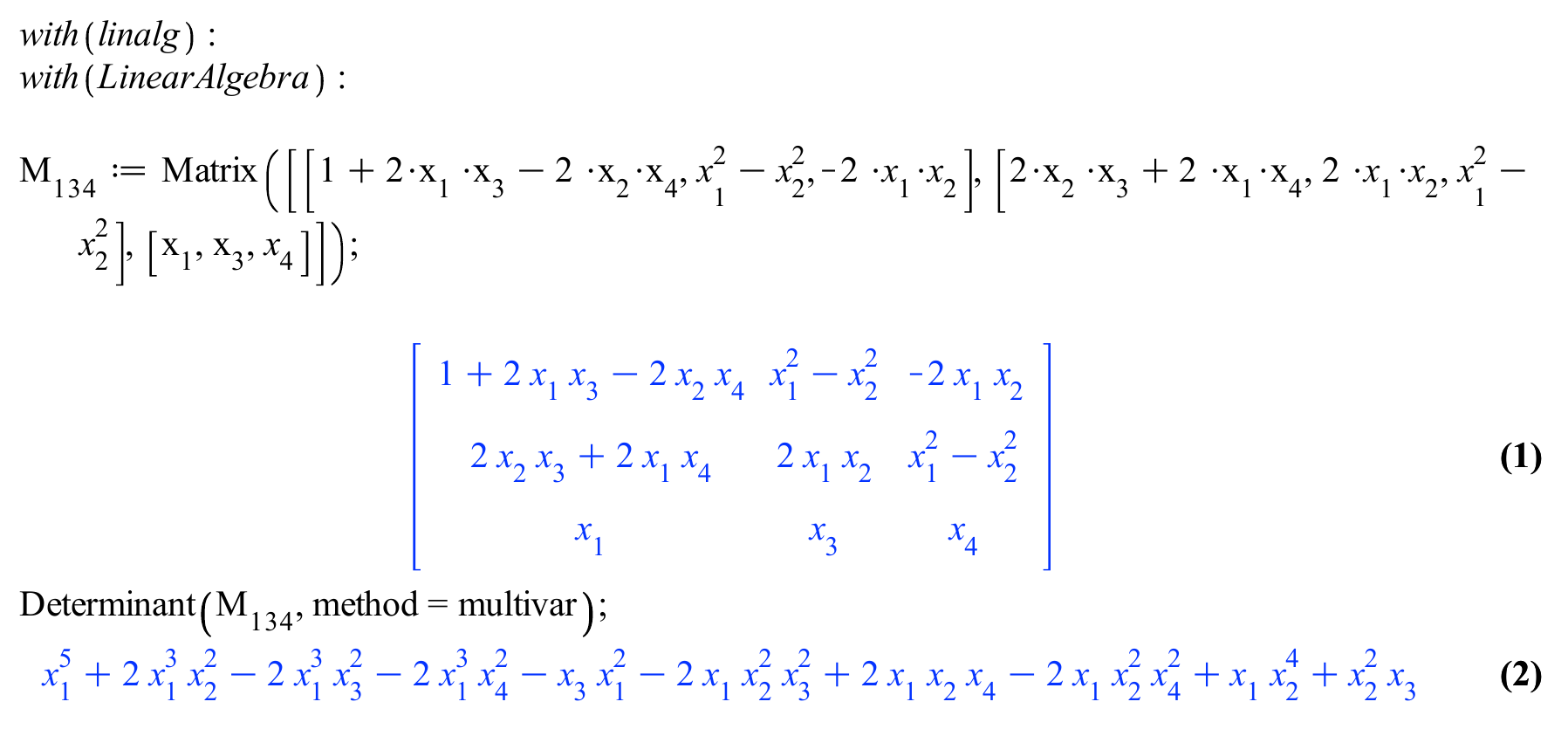}
\end{figure}

4) Calculate the determinant of the minor defined by the columns 2, 3 and 4:

\begin{figure}[h!]
\includegraphics[scale=0.6]{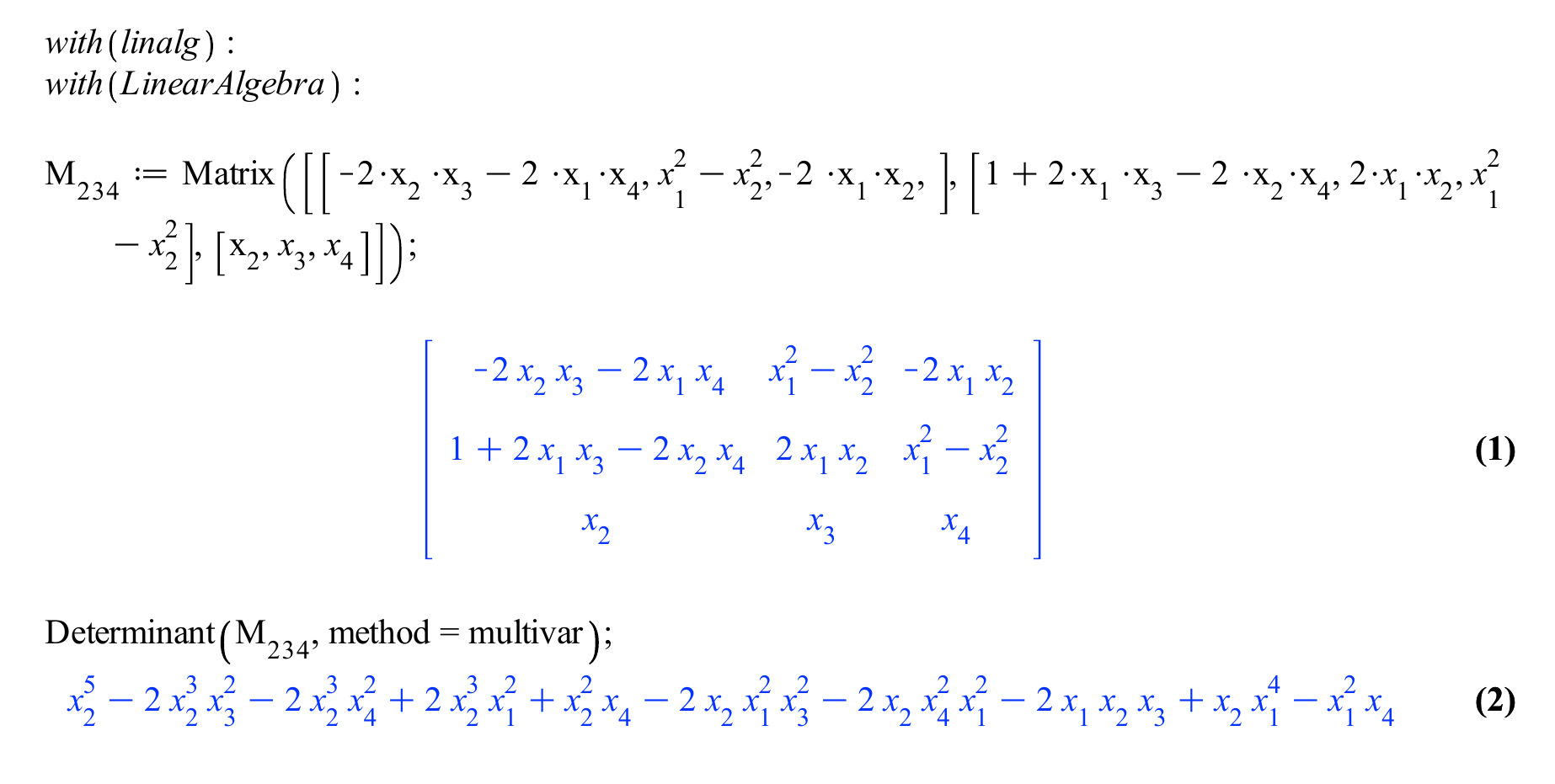}
\end{figure}

\pagebreak

B) Solve now the system of equations of the above determinants:

\begin{figure}[h!] 
\includegraphics[scale=0.6]{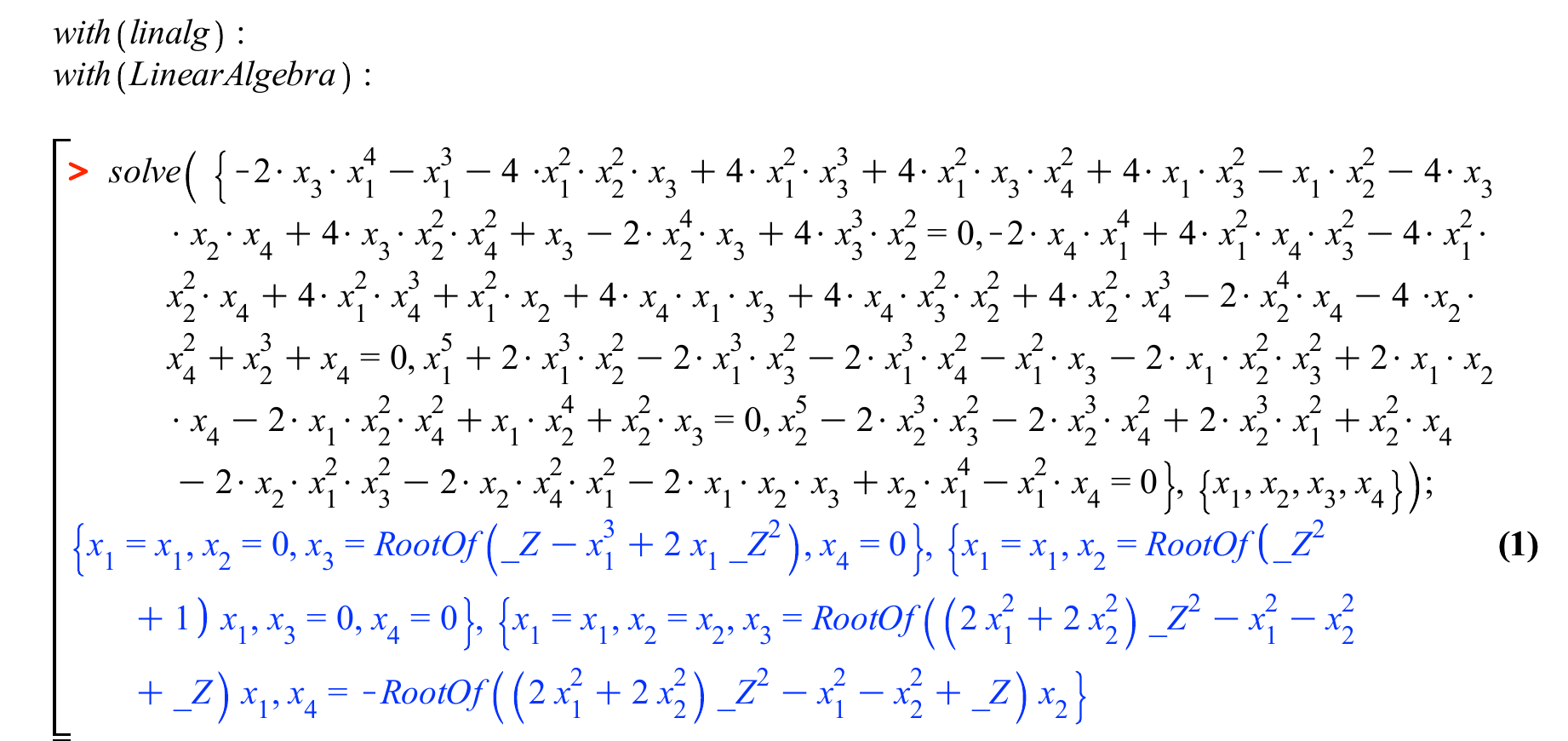}
\end{figure}

We conclude that  ${\mathcal{N}}_G = N_1 \cup N_2 \cup N_3,$
where 

\begin{figure}[h!] 
\includegraphics[scale=0.4]{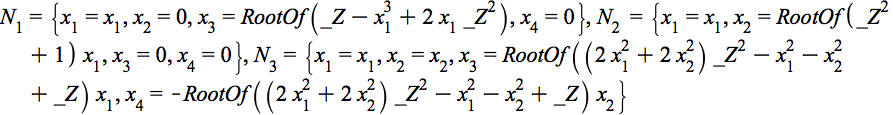}.
\end{figure}

C) In order to calculate ${\mathcal{V}}_G$, we have to calculate and draw $(G, \varphi) (N_i)$, for $i = 1, 2, 3$. 

1) Calculate and draw $(G, \varphi) (N_1)$: 
\begin{figure}[h!] 
\includegraphics[scale=0.6]{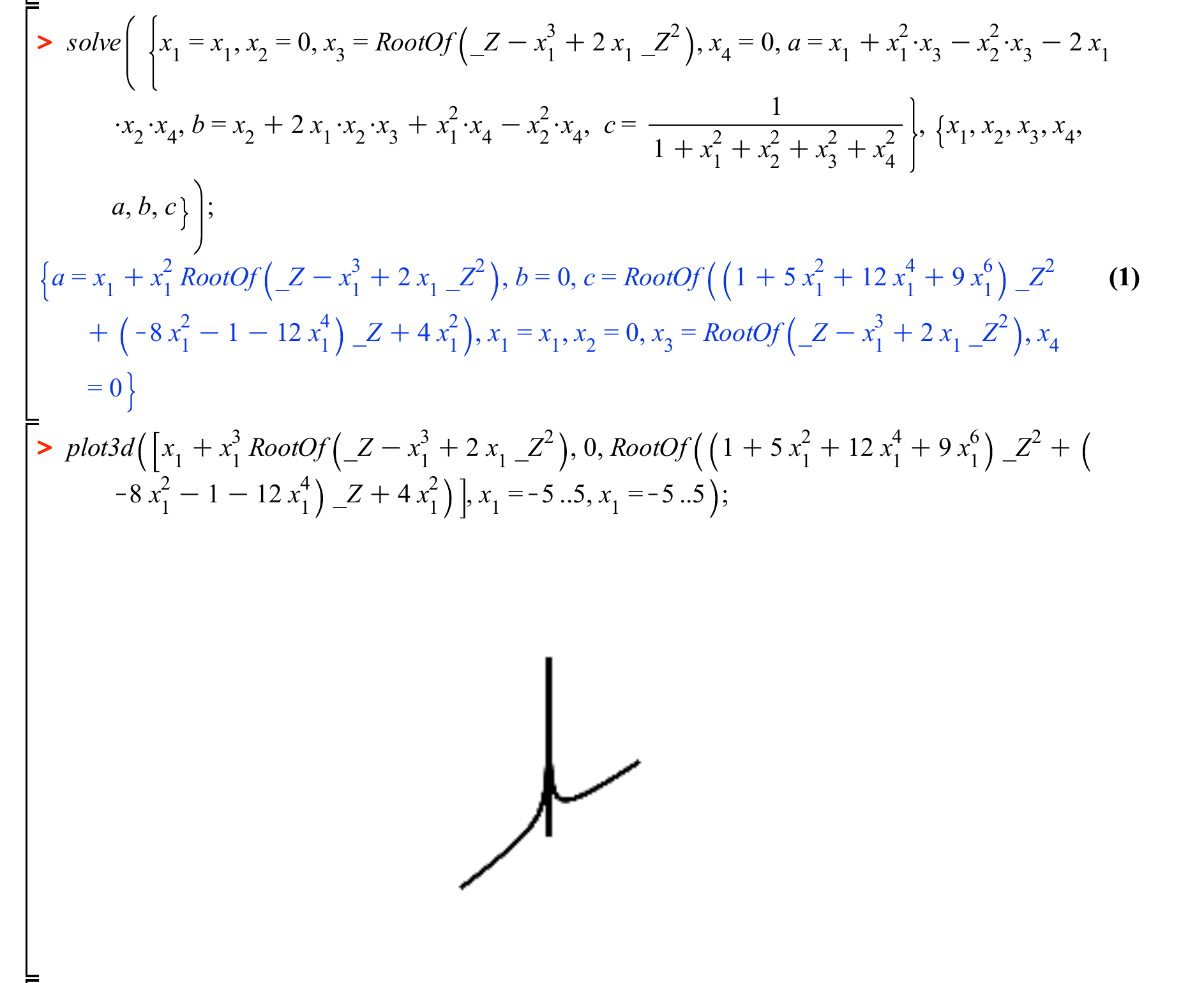}
\end{figure}

\pagebreak

2) Calculate and draw $(G, \varphi) (N_2)$: 

\begin{figure}[h!] 
\includegraphics[scale=0.6]{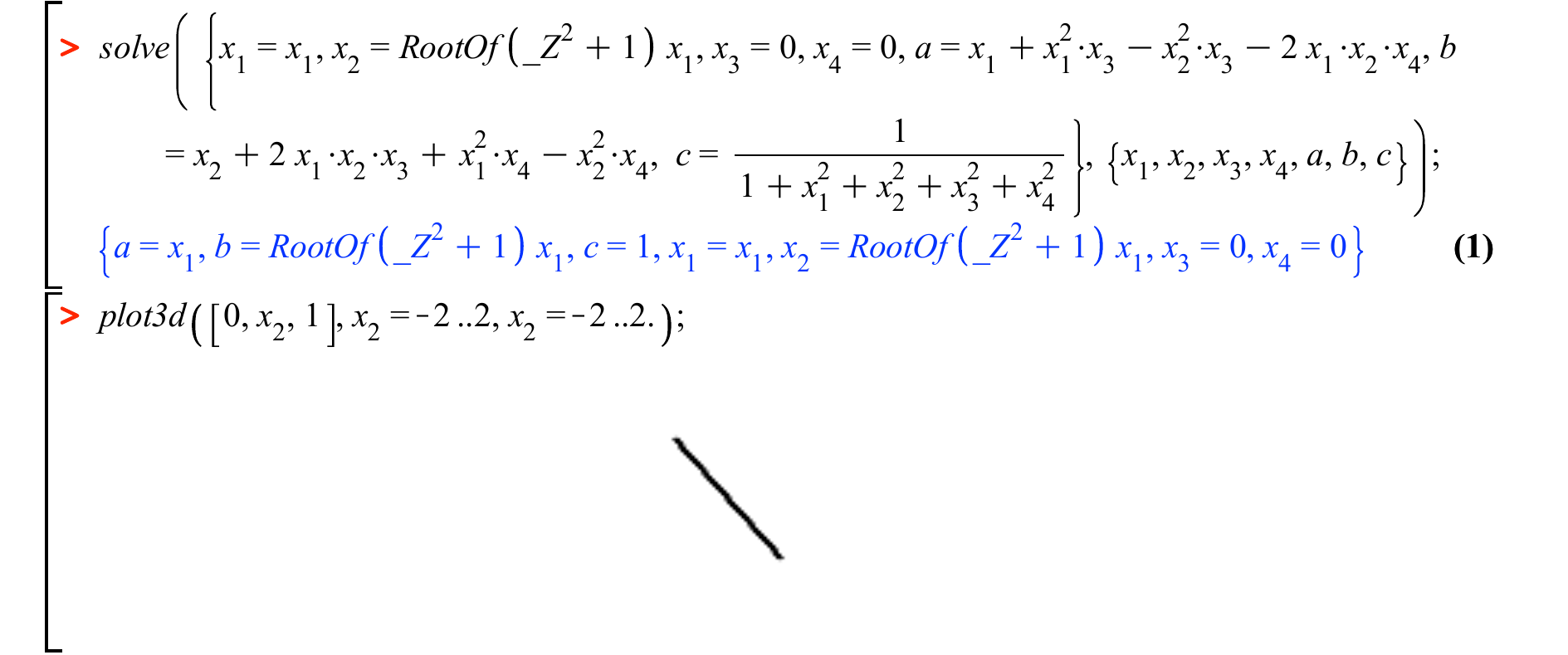}
\end{figure}

+ Calculate and draw $(G, \varphi) (N_3)$: 

\begin{figure}[h!] 
\includegraphics[scale=0.6]{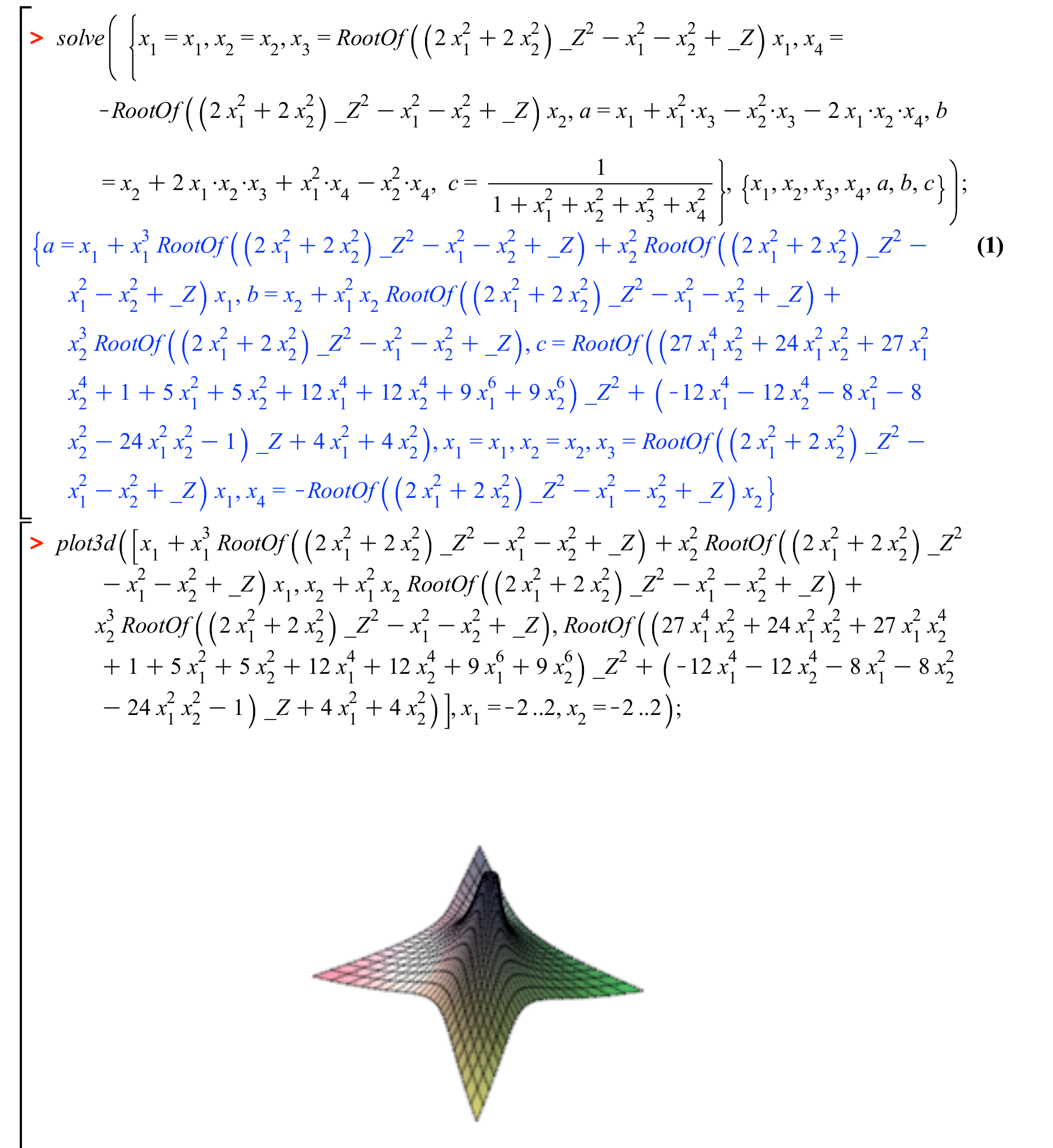}
\end{figure}

\pagebreak 

Since ${\mathcal{V}}_G$ is the closure of 
$\bigcup_{i = 1}^3 (G, \varphi) (N_i)$ then   ${\mathcal{V}}_G$ is connected and has a pure dimension, then $V_F$ is a cone:

\begin{figure}[h!]  
\includegraphics[scale=0.6]{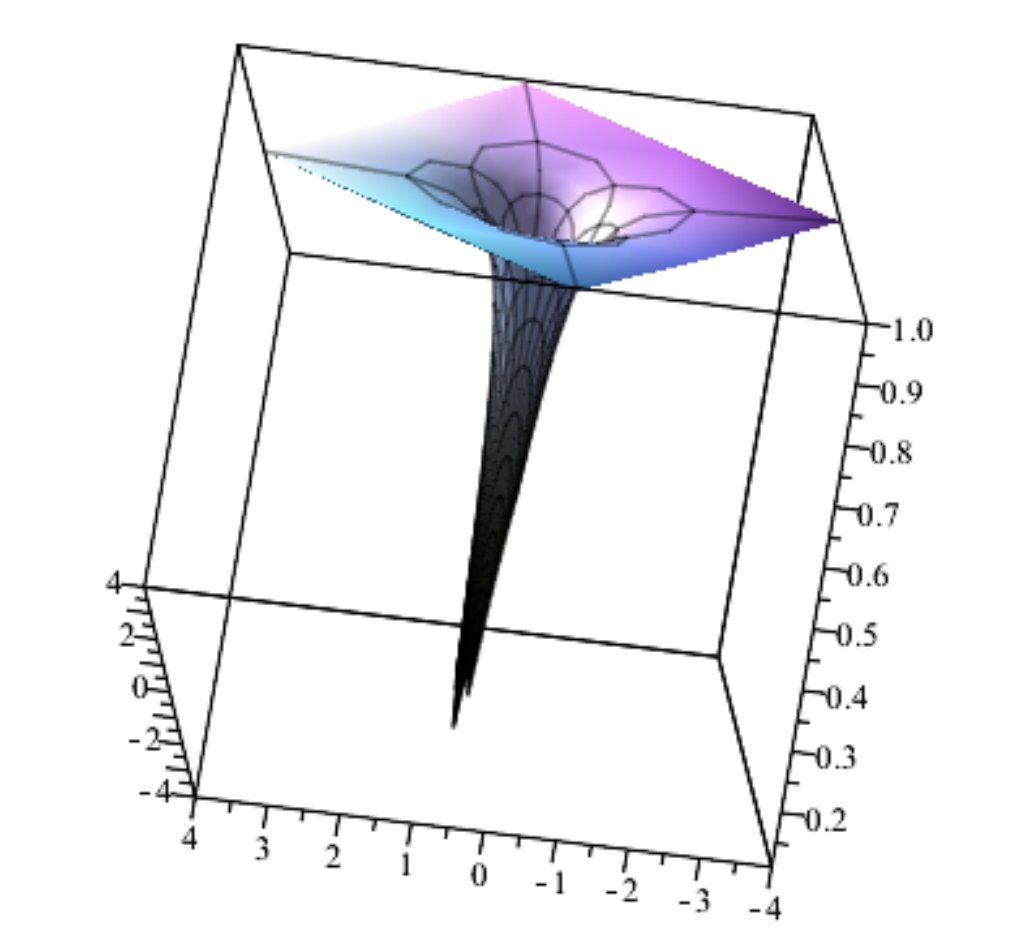} \label{cone}
\end{figure}
}

\end{example}

\begin{example}
{\rm If we take {\it the suspension} of the Broughton's example 
$$G: \C^3 \to \C^2, \quad \quad \quad  G(z, w, \eta) =( z + z^2w, \eta),$$
then, by the same calculations than the example \ref{ex}, the variety  ${\mathcal{V}}_G$ is a cone as in the example \ref{ex} but it has dimension 4, in the space $\R^6$.   We can check easily that the  intersection homology in dimension 2 of the variety ${\mathcal{V}}_G$ of this example is non trivial. We get an example to illustrate the corollary \ref{coThuyCidinha1}.
}
\end{example}

\begin{example}
{\rm If we take the Broughton example for $G: \C^3 \to \C$ such that $G(z, w, \eta) = z + z^2w,$  then similarly to the example  \ref{ex}, we get an example of varieties ${\mathcal{V}}_G$ for the case $G: \C^n \to \C^m$ where $ m \leq n-2$. This example illustrate the remark \ref{remark2}.
}
\end{example}

\bibliographystyle{plain}

\end{document}